\documentclass[preprint]{autart}    % Enable this line and disable the 
\newcommand{\EQ}{\begin{eqnarray}}
\newcommand{\EN}{\end{eqnarray}}
\newcommand{\EQQ}{\begin{eqnarray*}}
\newcommand{\ENN}{\end{eqnarray*}}
\newcommand{\nnum}{\nonumber}
\renewcommand{\theequation}{\thesection.\arabic{equation}}

\newcommand{\bremark}{\begin{remark} \begin{rm} }
\newcommand{\eremark}{ \end{rm}  \qed\\
\end{remark} }
\newcommand{\btheorem}{\begin{thm} \begin{rm} }
\newcommand{\etheorem}{ \end{rm}  \qed\\
\end{thm} }
\newcommand{\blemma}{\begin{lemma} \begin{rm} }
\newcommand{\elemma}{ \end{rm}  \qed\\
\end{lemma} }
\newcommand{\bcorollary}{\begin{corollary} \begin{rm} }
\newcommand{\ecorollary}{ \end{rm}  \qed\\
\end{corollary} }
\newcommand{\halmos}	{${ }^{ }$\hfill\rule{2mm}{2mm}}

\newcommand{\eps}				{{\epsilon}}

\newcommand{\R}		{{\mathds{R}}}
\newcommand{\er}[1]	{{(\ref{#1})}}

\usepackage{textcomp}
\usepackage{graphicx}          % Include this line if your 
   \usepackage{cite}                        % document contains figures,
\usepackage{amsmath}   % or this line, depending on which
\usepackage{amssymb}                              % you prefer.
\usepackage{dsfont}
\usepackage{eufrak}  

\usepackage{eufrak}      
\usepackage{dsfont}
\usepackage{cite}       
\usepackage{graphics} 
\usepackage{epsfig} 
\usepackage{times} 
\usepackage{pifont}
\usepackage{enumerate}
\usepackage{mathrsfs}

\begin{document}
\begin{frontmatter}
%\runtitle{Insert a suggested running title}  % Running title for regular 
                                              % papers but only if the title  
                                              % is over 5 words. Running title 
                                              % is not shown in output.

%\title{In Catilinam IV\thanksref{footnoteinfo}} % Title, preferably not more 
\title{A local characterization of Lyapunov functions on Riemannian manifolds $^{\star}$}                                                % than 10 words.

\thanks[footnoteinfo]{This work was supported by the Australian Research Council Discovery Project DP120101144.}

\author[First]{Farzin Taringoo}\ead{ftaringoo@unimelb.edu.au},
\author[First]{Peter M. Dower}\ead{pdower@unimelb.edu.au},    % Add the 
\author[First]{Dragan Ne\v{s}i\'{c}}\ead{dnesic@unimelb.edu.au}, 
\author[First]{Ying Tan}\ead{yingt@unimelb.edu.au}
             % e-mail address 
%\author[Baiae]{Publius Maro Vergilius}\ead{vergilius@culture.ir}  % (ead) as shown

\address[First]{Electrical and Electronic Engineering Department, The University of
Melbourne, Victoria, 3010, Australia} 
%\address[Second]{Pizza Delivery, Melbourne, 3010, Australia}   % Please supply                                              
%\address[Rome]{Senate House, Rome}             % full addresses
%\address[Baiae]{The White House, Baiae}        % here.

                             % keyword list or with the 
                                          % help of the Automatica 
                                          % keyword wizard

\begin{abstract}                          % Abstract of not more than 200 words.
This  paper proposes  converse Lyapunov theorems for  nonlinear dynamical systems defined on smooth connected Riemannian manifolds and characterizes properties of  Lyapunov functions with respect to the Riemannian distance function. We extend classical Lyapunov converse theorems for dynamical systems in $\mathds{R}^{n}$ to dynamical systems evolving on Riemannian manifolds. This is performed by restricting our analysis to the so called normal neighborhoods of equilibriums on Riemannian manifolds. By employing the derived properties of Lyapunov functions, we obtain the stability of perturbed dynamical systems on Riemannian manifolds.\end{abstract}

\end{frontmatter}

\section{Introduction}

Many systems include dynamics that naturally evolve on Riemannian manifolds, see for example \cite{Lewis, Bloch, Arnold, mar1, Sastry}, with their analysis requires the application of differential geometric tools. Examples of such systems can be found in many mechanical settings, see \cite{Lewis, Bloch, Arnold}.

Stability theory is an important topic in control theory. This theory addresses the stability of trajectories of dynamical systems as solutions of differential equations or differential inclusions, see \cite{Kha, Sastry, shu}. Lyapunov stability theory is the core mathematical tool for analyzing and characterizing the stability of equilibria. Stability in the sense of Lyapunov has been extensively analyzed in the literature, see for example \cite{Kha, las,lyap}.

Traditionally, the development of stability theory has focused on dynamical systems evolving on Euclidean spaces.
There, the application of the attendant vector space properties of Euclidean spaces leads to significant
simplifications in the analysis, due to coordinate transformations available to shift equilibria to the origin.
Consequently, stability analysis can be reduced to the analysis of a equilibria located at the origin. However, many dynamical systems defined on manifolds do not necessarily possess the vector space properties of Euclidean spaces. Consequently, a generalization of the traditional framework of the stability theory is inevitable. 

Stability analysis for systems evolving on manifolds is a new area of research, see for example \cite{mou,mal}.
Recent results concerning the existence and properties of Lyapunov functions are documented in \cite{hur,pat,gre,shu,kel,lin,haf,wil,sou,con,his}. In particular, in \cite{sou,con}, the existence of complete Lyapunov functions for dynamical systems on compact metric spaces is derived. In general, Riemannian manifolds can be considered as metric spaces by employing the notion of Riemannian distance function, see \cite{Lee3}.

In this paper, we present several converse Lyapunov theorems for dynamical systems evolving on Riemannian manifolds  and prove some local properties of such Lyapunov functions. To this end, we define Lyapunov stability of dynamical systems on Riemannian manifolds based on the Riemannian distance function.
We employ the notion of geodesics on Riemannian manifolds and apply the stability results for dynamical systems on $\mathds{R}^{n}$ to obtain the existence of Lyapunov functions for dynamical systems defined on Riemannian manifolds.
Using a version of stability theory for systems evolving on Riemannian manifolds, see \cite{forni,Ang, Lewis},  the stability results for dynamical systems evolving on Euclidean spaces \cite{Kha, nes, Hahn} are extended to those evolving on Riemannian manifolds. We introduce a lift operator to convert the dynamical equations on a Riemannian manifold to a dynamical system on the tangent space of an equilibrium, see \cite{Lee3, jost, Pet}, and invoke some of the standard results of the stability theory presented in \cite{Kha, Hahn}. It is shown that in a \textit{normal} neighborhood \cite{Lee3} of an equilibrium of a dynamical system, the constructed Lyapunov functions satisfy certain properties which can be used to analyze the stability and robustness of the underlying dynamical system. These results are extended and applied to
study perturbed dynamic systems. Geometric features of the normal neighborhoods, such as existence of unique length minimizing geodesics and their local representations enable us to closely relate the stability results obtained for dynamical systems in $\mathds{R}^{n}$ to those defined on Riemannian manifolds.

In terms of exposition,  Section \ref{s2}  presents some  mathematical preliminaries needed for the subsequent analysis. Section \ref{s3} presents the main results for the existence of Lyapunov functions for dynamical systems evolving on Riemannian manifolds.  These results are employed in Section \ref{s4}  to derive the stability of  perturbed dynamical systems on Riemannian manifolds. The paper concludes with some closing remarks in Section \ref{s5}.
\section{Preliminaries}\label{s2}In this section we provide the differential geometric material which is necessary for the analysis presented in the rest of the paper.  We define some of the frequently used symbols of this paper in Table 1.\\
\begin{table}[ht]
%\label{table1}
\caption{Symbols and Their Descriptions} % title of Table
\centering % used for centering table
\begin{tabular}{c c} % centered columns (4 columns)
\hline\hline %inserts double horizontal lines
Symbol& Description  \\ [0.5ex] % inserts table
%heading
\hline % inserts single horizontal line
$M$ & Riemannian manifold \\ 
$\mathfrak{X}(M)$ & space of smooth time invariant\\ &vector fields on $M$\\ 
$\mathfrak{X}(M\times \mathds{R})$ & space of smooth time varying\\ &vector fields on $M$\\ 
% inserting body of the table
$T_{x}M$ & tangent space at $x\in M$  \\
$T^{*}_{x}M$ & cotangent space at $x\in M$  \\
$TM$ & tangent bundle of $M$\\
$T^{*}M$ & cotangent bundle of $M$\\
$\frac{\partial}{\partial x_{i}}$ & basis tangent vectors at $x\in M$\\
$dx_{i}$ & basis cotangent vectors at $x\in M$\\
$f(x,t)$ & time-varying vector fields on $M$  \\
$||\cdot||_{g}$ & Riemannian norm\\
$||\cdot||_{e}$ & Euclidean norm\\
$||\cdot||$ & induced norm\\
$g(\cdot,\cdot)$ & Riemannian metric on $M$ \\ 
$d(\cdot,\cdot)$ & Riemannian distance on $M$\\
$\Phi_{f}$ & flow associated with $f$\\
$TF$ & pushforward of $F$\\
%$T^{*}\Phi_{f}$ & pull-back of $\Phi_{f}$ \\
$T_{x}F$ & pushforward of $F$ at $x$\\
%$T_{x}^{*}\Phi_{f}$ & pull-back of $\Phi_{f}$ x$\\
$\R_{>0}$& $(0,\infty)$\\
$\R_{\geq 0}$& $[0,\infty)$\\
$C^{\infty}(M)$ & space of smooth functions on $M$\\
$\simeq$ & isomorphism\\
$B(x,r)$ & metric ball centered at $x$ with radius $r$\\
$B_{r}(0)$ & Ball with radius $r$ in tangent spaces\\
\hline %inserts single line
\end{tabular}
\label{table:nonlin} % is used to refer this table in the text
\end{table}
\newtheorem{theorem}{Theorem}
\newtheorem{definition}{Definition}
\begin{definition}
\label{cc}
Let $M$ be a an $n$ dimensional manifold. A \textit{coordinate chart} on $M$ is $(U,\phi)$, where $U$ is an open set in $M$ and $\phi$ is a homomorphism from $U$
to $\phi(U)\subset \mathds{R}^{n}$, see \cite{Lee2}.
\halmos\end{definition}
\subsection{Riemannian manifolds}
\begin{definition}[see \cite{Lee2}, Chapter 3]
\label{def1}
A Riemannian manifold $(M,g)$ is a differentiable manifold $M$ together with a Riemannian metric $g$, where $g$ is defined for each $x\in M$ via an inner product $g_x:T_x M\times T_x M\rightarrow\R$ on the tangent space $T_x M$ (to $M$ at $x$) such that the function defined by $x\mapsto g_x(X(x),Y(x))$ is smooth for any vector fields $X,Y\in\mathfrak{X}(M)$. In addition,
\begin{enumerate}[(i)]
\item $(M,g)$ is $n$ dimensional if $M$ is $n$ dimensional;
\item $(M,g)$ is connected if for any $x,y\in M$, there exists a piecewise smooth curve that connects $x$ to $y$. 
\end{enumerate}
\halmos
\end{definition}
Note that in the special case where $M\doteq\mathds{R}^{n}$, the Riemannian metric $g$ is defined everywhere by $g_{x}=\sum^{n}_{i,j=1}g_{ij}(x)dx_{i}\otimes dx_{j}$, where $\otimes $ is the tensor product on $T^{*}_{x}M\times T^{*}_{x}M $, see \cite{Lee2}.

As formalized in Definition \ref{def1}, connected Riemannian manifolds possess the property that any pair of points $x,y\in M$ can be connected via a path $\gamma\in\mathscr{P}(x,y)$, where
\EQ
	\mathscr{P}(x,y)
	\doteq \left\{ \gamma:[a,b]\rightarrow M \, \left| \, \begin{array}{c} 
			\gamma \mbox{ piecewise smooth,}
			% \\
			% \text{with}
			\\
			\gamma(a) = x\,, \ \gamma(b) = y
		\end{array} \right. \right\}
	\label{eq:paths}
\EN

\begin{theorem}[ \hspace{-.1cm}\cite{Lee3}, P. 94]
\label{t1}
\label{thm:path}
Suppose $(M,g)$ is an $n$ dimensional connected Riemannian manifold. Then, for any $x,y\in M$, there exists a piecewise smooth path $\gamma\in\mathscr{P}(x,y)$ that connects $x$ to $y$.
\halmos
\end{theorem}

The existence of connecting paths (via Theorem \ref{thm:path}) between pairs of elements of an $n$ dimensional connected Riemannian manifold $(M,g)$ facilitates the definition of a corresponding Riemannian distance. In particular, the Riemannian distance $d:M\times M\rightarrow\R$ is defined by the infimal path length between any two elements of $M$, with
\EQ\label{len}
	d(x,y)
	& \doteq \inf_{\gamma\in\mathscr{P}(x,y)} \int_a^b
						\sqrt{g_{\gamma(t)} (\dot\gamma(t),\, \dot\gamma(t))}\, dt\,.
	\label{eq:distance-via-paths}
\EN
Note that in the special case where $M\doteq\R^n$, the Riemannian distance \er{eq:distance-via-paths} simplifies to $d(x,y) = \|x - y\|_{e}$.

Using the definition of Riemannian distance $d$ of \er{eq:distance-via-paths}, $(M,d)$ defines a metric space as formalized by the following theorem.
\begin{theorem}[ \hspace{-.1cm}\cite{Lee3}, P. 94]
\label{t2}
Any $n$ dimensional connected Riemannian manifold $(M,g)$ defines a metric space $(M,d)$ via the Riemannian distance $d$ of \er{eq:distance-via-paths}. Furthermore, the induced topology of $(M,d)$ is the same as the manifold topology of $(M,g)$. 
\halmos
\end{theorem}
Next, the crucial \textit{pushforward} operator is introduced.
\begin{definition}
\label{d2}
For a given smooth mapping $F:M\rightarrow N$ from manifold $M$ to manifold $N$ the pushforward $TF$ is defined as a generalization of the Jacobian of smooth maps between Euclidean spaces as follows:
\EQ TF:TM\rightarrow TN,\EN
where 
\EQ T_{x}F:T_{x}M\rightarrow T_{F(x)}N,\EN
and 
\EQ T_{x}F(X_{x})\circ h=X_{x}(h\circ F),\hspace{.2cm}X_{x}\in T_{x}M, h\in C^{\infty}(N).\nnum\\\EN
\halmos\end{definition}
\newtheorem{lemma}{Lemma}
%\begin{lemma}(\cite{Lee3}, Page 50)

%\label{l1}
%Consider $\gamma:(-\epsilon,\epsilon)\rightarrow M$ such that $\gamma(0)=p\in M$ and $\dot{\gamma}(0)=X_{p}$. If two vector fields $Y$ and $\tilde{Y}$ agree along $\gamma$, then 
%\EQ\nabla_{X_{p}}Y|_{p}= \nabla_{X_{p}}\tilde{Y}|_{p}.\EN\hspace*{7.8cm}\halmos
%\end{lemma}

%Since $\nabla$ is torsion free, i.e. $\nabla_{X}Y-\nabla_{Y}X=[X,Y]$ and $[\frac{\partial}{\partial \epsilon},\frac{\partial}{\partial \tau}]=0$, we have
%\EQ \nabla_{\frac{\partial}{\partial \epsilon}}\partial_{\tau}\Gamma(\epsilon,\tau)=\nabla_{\frac{\partial}{\partial \tau}}\partial_{\epsilon}\Gamma(\epsilon,\tau).\EN
%The property above will be used to extend the standard averaging techniques for dynamical systems defined on Riemannian manifolds.
\subsection{Dynamical systems on Riemannian manifolds}
This paper focuses on dynamical systems governed by differential equations on a connected $n$ dimensional Riemannian manifold $M$. Locally these differential equations are defined by (see \cite{Lee2})
\EQ \label{peter}&&\hspace{-0cm}\dot{x}(t)=f(x(t),t),\hspace{.2cm} f\in \mathfrak{X}(M\times \mathds{R}),\nnum\\&&\hspace{-0cm}\hspace{.2cm} x(0)=x_{0}\in M, t\in[t_{0},t_{f}].\EN
The time dependent flow associated with a differentiable time dependent vector field $f$ is a map $\Phi_{f}$ satisfying
\EQ \label{flow} &&\Phi_{f}:[t_0, t_{f}]\times [t_{0}, t_{f}]\times M\rightarrow M, \nnum\\&& (s_{0},s_{f},x)\mapsto \Phi_{f}(s_{f},s_{0},x)\in M,\EN
and
\EQ \left.\frac{d\Phi_{f}(s,s_{0},x)}{ds}\right|_{s=t}=f(\Phi_{f}(t,s_{0},x),t).\EN 
One may show, for a smooth vector field $f$, the integral flow $\Phi_{f}(s,t_{0},\cdot):M\rightarrow M$ is a local diffeomorphism , see \cite{Lee2}.
  Here we assume that the vector field $f$ is smooth and \textit{complete}, i.e. $\Phi_{f}$ exists for all $t\in (t_{0},\infty)$.
  \subsection{Geodesic Curves}
   Geodesics are defined  \cite{jost} as length minimizing curves on Riemannian manifolds which satisfy
   \EQ \nabla_{\dot{\gamma}(t)}\dot{\gamma}(t)=0,\EN
   where $\gamma(\cdot)$ is a geodesic curve on $(M,g)$ and $\nabla$ is the \textit{Levi-Civita} connection on $M$, see \cite{Lee3}.
    The solution of the Euler-Lagrange variational problem associated with the length minimizing problem shows that  all the geodesics on an $n$ dimensional Riemannian manifold  $(M,g)$  must satisfy the following system of ordinary differential equations:
\EQ \label{geo}\ddot{\gamma}_{i}(s)+\sum^{n}_{j,k=1}\Gamma^{i}_{j,k}\dot{\gamma}_{j}(s)\dot{\gamma}_{k}(s)=0,\quad i=1,...,n,\EN
where
\EQ\label{cris} \Gamma^{i}_{j,k}=\frac{1}{2}\sum^{n}_{l=1}g^{il}(g_{jl,k}+g_{kl,j}-g_{jk,l}),\quad g_{jl,k}=\frac{\partial g_{jl}}{\partial x_{k}},\nnum\\\EN
where all the indexes  $i,j,k,l$ run from $1$ up to $n=dim(M)$ and $[g^{ij}]\doteq[g_{ij}]^{-1}$. Note that $g_{ij}$ is the $(i,j)$ entity of the metric $g$.
 \begin{definition}[ \hspace{-.1cm}\cite{Lee3}, p. 72]
 The restricted exponential map is defined by 
 \EQ \exp_{x}:T_{x}M\rightarrow M,\hspace{.2cm}\exp_{x}(v)=\gamma_{v}(1), v\in T_{x}M,\EN
 where $\gamma_{v}(1)$ is the unique maximal geodesic \cite{Lee3}, P. 59, initiating from $x$ with the velocity $v$ up to one.
 \halmos\end{definition}
 Throughout, \textit{restricted exponential maps} are referred to as \textit{exponential maps}. 
 An open ball of radius $\delta>0$ and centered at $0\in T_x M$ in the tangent space at $x$ is denoted 
by $B_{\delta}(0)\doteq\{v\in T_{x}M\hspace{.1cm}|\hspace{.2cm}||v||_{g}<\delta\}$. Similarly, the corresponding closed ball is denoted by $\overline{B}_\delta(0)$. Using the local  diffeomorphic property of exponential maps, the corresponding geodesic ball centered at $x$ is deﬁned as follows. 
 
 \begin{lemma}[ \hspace{-.1cm}\cite{Lee3}, Lemma 5.10]
 \label{eun}
 For any $x\in M$, there exists a neighborhood $B_{\delta}(0)$ in $T_{x}M$ on which $\exp_{x}$ is a diffeomorphism onto $\exp_{x}(B_{\delta}(0))\subset M$. 
 \halmos\end{lemma}
 \begin{definition}[\hspace{-.01cm}\cite{Lee3}]
 In a neighborhood of  $x\in M$, where $\exp_{x}$ is a local diffeomorphism (this neighborhood always exists by Lemma \ref{eun}), a geodesic ball of radius $\delta>0$ is denoted by $\exp_{x}(B_{\delta}(0))\subset M$. The corresponding closed geodesic ball is denoted by $\exp_{x}(\overline{B}_{\delta}(0))$. 
 \halmos\end{definition}
 
 \begin{definition}
 For a vector space $V$, a \textit{star-shaped neighborhood} of $0\in V$ is any open set $U$ such that if $u\in U$ then $\alpha u\in U, \alpha\in[0,1]$.
 \halmos\end{definition}
 \begin{definition}[ \hspace{-.1cm}\cite{Lee3}, p. 76]
 A normal neighborhood around $x\in M$ is any open neighborhood  of $x$ which is a diffeomorphic image of a star shaped neighborhood of $0\in T_{x}M$ under $\exp_{x}$ map.
 \halmos\end{definition}
 %\begin{definition}[ \hspace{-.1cm}\cite{Lee3}]
 %A uniformly normal neighborhood around $x\in M$ is any open neighborhood $\mathcal{U}_{x}$ of $x$ for which there exists $\delta>0$, such that $\mathcal{U}_{x}\subset \exp_{p}B_{\delta}(0)$ for all $p\in \mathcal{U}_{x}$.
 %\halmos\end{definition}

 %\begin{lemma}[\hspace{-.01cm}\cite{Lee3}]
 %\label{un}
 %For any $x\in M$ and any neighborhood $\mathcal{U}_{x}$ of $x$, there exists a uniformly normal neighborhood $\mathcal{V}_{x}$ such that $\mathcal{V}_{x}\subset \mathcal{U}_{x}$.
 %\halmos\end{lemma}
 \begin{definition}\label{inj}The injectivity radius of $M$ is  
 \EQ i(M)\doteq \inf_{x\in M}i(x),\EN
 where
\EQ&&  i(x) \doteq \sup\{ r\in\R_{\ge0}| \exp_x \mbox{is diffeomorphic onto}\nnum\\&& \exp_x (B_r(0))\}.\nnum\\&&\hspace{-1cm}\EN
\halmos
\end{definition}
\begin{definition}
The metric ball with respect to $d$ on $(M, g)$ is defined by
\EQ B(x,r)\doteq \{ y\in M\hspace{.1cm}|\hspace{.2cm}d(x,y)<r\}.\EN
\halmos\end{definition}
The following lemma reveals a relationship between normal neighborhoods and metric balls on $(M,g)$.
\begin{lemma}[ \hspace{-.1cm}\cite{Pet}, p. 122]
\label{lpp}
Given any $\eps\in\R_{>0}$ and $x\in M$, suppose that $\exp_x$ is a diffeomorphism  on $B_{\epsilon}(0)\subset T_{x}M$, and $B(x,r)\subset \exp_{x}B_{\epsilon}(0)$ for some $r\in \mathds{R}_{>0}$. Then 
\EQ \exp_{x}B_{r}(0)=B(x,r).\EN
\hspace*{7.5cm}\halmos\end{lemma}
We note that $B_{\epsilon}(0)$ is the metric ball of radius $\epsilon$ with respect to the Riemannian metric $g$ in $T_{x}M$.
%\begin{lemma}[\hspace{-.01cm}\cite{Klin}]
%\label{kl}
%The injectivity radius $i(x), x\in M$ is continuous with respect to $x$ and is bounded below for compact Riemannian manifolds.
%\halmos\end{lemma}
 
  %%%%%%%%%%%%%%%%%%%%%%%
%%%%%%%%%%%%%%%%%%%%%%%%%%
%%%%%%%%%%%%%%%%%%%%%%%%%%%%%
%%%%%%%%%%%%%%%%%%%%%%%%%%%%%%%%
%%%%%%%%%%%%%%%%%%%%%%%%%%%%\section{Hybrid systems}
%%%%%%%%%%%%%%%%%%%%%%%%%%%%%%%%%%%%
%%%%%%%%%%%%%%%%%%%%%%%%%%%%%%%%%%%

 %%%%%%%%%%%%%%%%%%%%%%%%%
 \section{Lyapunov Analysis on Riemannian Manifolds}
 \label{s3}
 We  extend the notion of stability to dynamical systems evolving on Riemannian manifolds. This problem has been addressed in \cite{mar1,Lewis,leca} in a geometric framework. The main motivation here is to characterize the local properties of Lyapunov functions based upon the Riemannian distance function. These properties will be of great importance in analyzing a range of dynamical systems evolving on manifolds. 
 
\begin{figure}
%\begin{figure}
  \vspace*{0cm}\begin{center}
\hspace*{-2cm}\includegraphics[scale=.25]{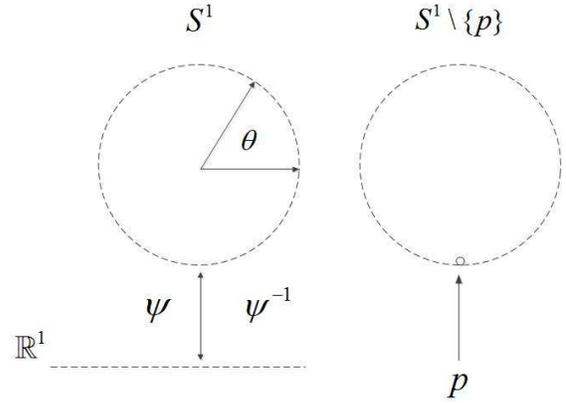}
 \caption{$S^{1}$ and $S^1\setminus\{p\}$}
     \label{ff12}
     \end{center}
 %  \end{figure}
  \end{figure}

 It is important to note that, depending on the geometry of the state space of a particular dynamical system, Riemannian distance might be significantly different than the Euclidean distance of embedded manifolds. As an example consider a unit circle $S^{1}\subset \mathds{R}^{2}$ in Figure \ref{ff12}. A \textit{local coordinate system} \cite{Lee4} for $S^{1}$ is given by the local homeomorphism $\psi:S^1\rightarrow\mathds{R}$ (see also Figure \ref{ff12}) defined by 
\EQ \theta \overset{\psi}{\mapsto} (\sin(\theta),\cos(\theta))\in \mathds{R}^{2}, \hspace{.2cm}\theta\in (0,2\pi)\subset \mathds{R}^{1}.\EN
In the case of the removal of a point $p$ from $S^{1}$, the Euclidean distance between points converging in $S^1\setminus\{p\}$ to $p\in S^1$ from either side converges to zero. However, at the same time, the Riemannian distance converges to $2\pi$ which is the largest distance on $S^{1}$ between any pair of points.

We generalize the stability notion for dynamical systems on Riemannian manifolds as follows. 

\begin{definition}
For the time-varying dynamical system $\dot{x}=f(x(t),t), \hspace{.2cm}f\in \mathfrak{X}(M\times \mathds{R})$,  $\bar{x}\in M$ is an equilibrium if 
\EQ \Phi_{f}(t,t_{0},\bar{x})=\bar{x},\hspace{.2cm}t\in[t_{0},\infty),\EN
where $\Phi_{f}$ is the integral flow of $f$ defined by (\ref{flow}).
\halmos\end{definition}
\begin{definition}[ \hspace{-.1cm}\cite{Lewis, Kha,forni,Ang}]
\label{les}
For the dynamical system $\dot{x}=f(x(t),t), \hspace{.2cm}f\in \mathfrak{X}(M\times \mathds{R})$, an equilibrium $\bar{x}\in M$ is\\ 

(i) \textit{uniformly Lyapunov stable} if  for any neighborhood $\mathcal{U}_{\bar{x}}$ of $\bar{x}\in M$ and any initial time $t_{0}\in \mathds{R}$, there exists a neighborhood $\mathcal{W}_{\bar{x}}$ of $\bar{x}$, such that 
\EQ \forall x_{0}\in \mathcal{W}_{\bar{x}}, \Phi_{f}(t,t_{0},x_{0})\in \mathcal{U}_{\bar{x}},\hspace{.2cm}\forall t\in[t_{0},\infty).\nnum\\\EN\\
(ii) uniformly locally asymptotically stable if it is Lyapunov stable and for any $t_0\in\mathds{R}$, there exists $\mathcal{U}_{\bar{x}}$ such that
\EQ &&\forall x_{0}\in \mathcal{U}_{\bar{x}},\hspace{.2cm}\lim_{t\rightarrow \infty}\Phi_{f}(t,t_{0},x_{0})=\bar{x},\hspace{.2cm}\mbox{i.e.}\nnum\\&&\lim_{t\rightarrow \infty}d(\Phi_{f}(t,t_{0},x_{0}),\bar{x})=0, \hspace{.2cm}t\in[t_{0},\infty). \EN\\
(iii) uniformly globally asymptotically stable if it is Lyapunov stable and for any $t_0\in\mathds{R}$, 
\EQ \forall x_{0}\in M, \hspace{.2cm}\lim_{t\rightarrow \infty}\Phi_{f}(t,t_{0},x_{0})=\bar{x}, \hspace{.2cm}t\in[t_{0},\infty). \EN
(iv) uniformly locally exponentially stable if it is locally asymptotically stable and for any $t_0\in\mathds{R}$, there exist $\mathcal{U}_{\bar{x}}$ and $K,\lambda\in \mathds{R}_{>0}$ such that
\EQ \label{es}&&\forall x_{0}\in \mathcal{U}_{\bar{x}}, d(\Phi_{f}(t,t_{0},x_{0}),\bar{x})\leq\nnum\\&& Kd(x_{0},\bar{x})\exp(-\lambda(t-t_{0})), t\in[t_{0},\infty).\EN
(v) globally exponentially stable if it is globally asymptotically stable and for any $t_0\in\mathds{R}$, there exist $K,\lambda\in \mathds{R}_{>0}$, such that,  
\EQ \label{es}&&\forall x_{0}\in M,d(\Phi_{f}(t,t_{0},x_{0}),\bar{x})\leq\nnum\\&& Kd(x_{0},\bar{x})\exp(-\lambda(t-t_{0})),t\in[t_{0},\infty).\EN\halmos
\end{definition}
We note that the convergence on $M$ is defined in the topology induced by $d$ which is the same as the original topology of $M$ by Theorem \ref{t2}.
\newtheorem{remark}{Remark}
\begin{definition}[\hspace{-.01cm}\cite{Lewis, Kha}]
A function $\chi:M\rightarrow \mathds{R}$ is locally positive definite (positive semi-definite) in a neighborhood of $\bar {x} \in M$  if $\chi(\bar{x})=0$ and there exists a neighborhood $\mathcal{U}_{\bar{x}}\subset M$ such that for all
$ x\in \mathcal{U}_{\bar{x}}\setminus\{\bar{x}\},\hspace{.2cm} 0<\chi(x) \hspace{.2cm}(\mbox{respectively}\hspace{.2cm} 0\leq \chi(x)).$\halmos
\end{definition}
Given a smooth function $\chi:M\rightarrow \mathds{R}$, the Lie derivative of $\chi$ along a time invariant vector field $f\in \mathfrak{X}(M)$ is defined by
\EQ \label{lie}\mathfrak{L}_{f}\chi\doteq d\chi(f),\EN
where $d\chi:TM\rightarrow \mathds{R}$ is the differential form of $\chi$. In any neighbourhood of $x\in M$, $d\chi$ is given locally by 
\EQ \label{pet2}d\chi=\sum^{n}_{i=1}\frac{\partial \chi}{\partial x_{i}}dx_{i}\in T^{*}_{x}M,\EN
where $n\doteq dim(M)$ and $T_x^* M$ is the cotangent space of $M$ at $x$, see \cite{Lee2}.

 \begin{remark}  For time-varying dynamical systems evolving on $M$, the Lie derivative of a smooth time-varying function $\chi:M\times \mathds{R}\mapsto \mathds{R}$ is defined by
\EQ\label{lie} \mathfrak{L}_{f(x,t)}\chi\doteq d\chi\left(\frac{\partial}{\partial t},f(x,t)\right),\EN
where 
\EQ d\chi=d_{x}\chi\oplus d_{t}\chi\in T^{*}_{x}M\oplus T^{*}_{t}\mathds{R},\nnum\\\EN
with $d_{x}\chi\in T^{*}_{x}M$ as per (\ref{pet2}), and $d_{t}\chi\in T^{*}_{t}\mathds{R}$.\end{remark}
\begin{definition}[\hspace{-.01cm}\cite{Kha, Lewis, mar1}](Lyapunov Candidate Functions)
\label{kir}
A smooth function $v:M\times \mathds{R}\rightarrow \mathds{R}$ is a Lyapunov function for the time-variant vector field $f\in \mathfrak{X}(M\times \mathds{R})$ if $v$ is locally positive definite in a neighborhood of an equilibrium $\bar{x}$ for $t\in [t_{0},\infty)$ and $\mathfrak{L}_{f}v$ is locally negative semi-definite in a neighborhood of $\bar{x}$.\halmos
\end{definition}
%\begin{remark}
%The Lyapunov function $v$ for the time-varying dynamical system $\dot{x}=f(x,t)$ satisfies the conditions in Definition \ref{kir} together with $v(\bar{x},t)=0,\forall t\in[t_{0},\infty)$, where $\bar{x}$ is an equilibrium of $f$.
%\end{remark}
%\begin{definition}
%The sublevel set $\mathcal{N}_{b}$ of a positive semidefinite function $v:M\rightarrow\mathds{R}$ is defined as $\mathcal{N}_{b}\doteq\{x\in M, \hspace{.2cm} v(x)\leq b\}$. By $\mathcal{N}_{b}(\bar{x})$ we denote the connected sublevel set of $M$ containing $\bar{x}\in M$. 
%\halmos\end{definition}%\begin{definition}
%\label{les}
%A smooth dynamical system $\dot{x}(t)=f(x,t),\hspace{.2cm}x(t)\in M$, is locally exponentially  stable with respect to the equilibrium $\bar{x}\in M$ if there exists $\mathcal{U}_{\bar{x}}\subset M$ such that 

%\begin{lemma}[ \hspace{-.1cm}\cite{Lewis}]
%\label{lc}
%Let $\bar{x}\in M$ be an equilibrium of $\dot{x}=f(x),\hspace{.2cm}x(t)\in M$ and 
%$v$ be a Lyapunov function on a neighborhood of $\bar{x}$. Then, for any neighborhood $\mathcal{U}_{\bar{x}}$ of $\bar{x}$, there exists $b\in\mathds{R}_{> 0}$ such that $\mathcal{N}_{b}(\bar{x})$ is compact, $\bar{x}\in int(\mathcal{N}_{b}(\bar{x}))$ and $\mathcal{N}_{b}(\bar{x})\subset \mathcal{U}_{\bar{x}}$, where $\mbox{int}(A)$ is the interior set  of $A$.\halmos  
%\end{lemma}
\begin{definition}
The time-variant sublevel set $\mathcal{N}_{b,t}$ of a positive semidefinite function $v:M\times \mathds{R}\rightarrow\mathds{R}$ is defined as $\mathcal{N}_{b,t}\doteq\{x\in M, \hspace{.2cm} v(x,t)\leq b\}$. By $\mathcal{N}_{b,t}(\bar{x})$ we denote a connected sublevel set of $M$ containing $\bar{x}\in M$. 
\halmos\end{definition}
The following lemma shows that there exists a connected compact sublevel set of an equilibrium point of a dynamical system on a Riemannian manifold.
\begin{lemma}
\label{lc}
Let $\bar x\in M$ and $v:M\times \mathds{R}\rightarrow\mathds{R}$ denote an equilibrium and a Lyapunov function 
respectively for system (\ref{peter}). Then, for any neighborhood $\mathcal{U}_{\bar{x}}$ of $\bar{x}$ and any $t\in \mathds{R}$, there exists $b\in\mathds{R}_{> 0}$, such that $\mathcal{N}_{b,t}(\bar{x})$ is compact, $\bar{x}\in int(\mathcal{N}_{b,t}(\bar{x}))$ and $\mathcal{N}_{b,t}(\bar{x})\subset \mathcal{U}_{\bar{x}}$, where $int(\cdot)$ gives the interior of a set.
\end{lemma}
\begin{pf}
The proof is based on the proof given in \cite{Lewis}, Lemma 6.12. In this case we fix time  $t\in\mathds{R}$ and consider $v(\cdot,t):M\rightarrow \mathds{R}$ as a smooth time-invariant function. In this case,  we apply the results of Lemma 6.12 in \cite{Lewis} to complete the proof.
\halmos\end{pf}
To analyze the behavior of dynamical systems on manifolds we employ the notion of \textit{comparison functions} defined in \cite{Kha}.
\begin{definition}[ \hspace{-.1cm}\cite{Kha}]
A continuous function $\alpha:[0,b)\rightarrow \mathds{R}_{\geq 0}$ is of class $\mathcal{K}$ if it is strictly increasing and $\alpha(0)=0$, and of class $\mathcal{K}_{\infty}$ if $b=\infty$ and  $\lim_{r\rightarrow\infty}\alpha(r)=\infty$.
\halmos\end{definition} 

\begin{definition}[ \hspace{-.1cm}\cite{Kha}]
A continuous function $\beta:[0,b)\times \mathds{R}_{\geq 0}\rightarrow \mathds{R}_{\geq 0}$ is of class $\mathcal{K}\mathcal{L}$ if for each fixed $s$, $\beta(\cdot,s)\in \mathcal{K}$ and for each fixed $r\in[0,b)$, $\beta(r,\cdot)$ is decreasing with $\lim_{s\rightarrow \infty}\beta(r,s)=0$. \halmos 
\end{definition} 

The following theorem provides $\mathcal{K}$ and $\mathcal{K}\mathcal{L}$ comparison function bounds for trajectories of uniformly stable dynamical systems evolving on Riemannian manifolds. 
\begin{theorem}
\label{lyp}
Any time-varying dynamical system of the form (\ref{peter}), evolving on a connected $n$ dimensional 
Riemannian manifold $(M,g)$, satisﬁes the following properties:
\begin{itemize}
 \item If an equilibrium $\bar{x}\in M$ is uniformly Lyapunov stable, then there exists a class $\mathcal{K}$ function $\alpha$ and a neighborhood $\mathcal{N}_{\bar{x}}$, such that 
\EQ &&d(\Phi_{f}(t,t_{0},x_{0}),\bar{x})\leq \alpha(d(x_{0},\bar{x})),\hspace{.2cm}\nnum\\&& x_{0}\in \mathcal{N}_{\bar{x}}, t\in[t_{0},\infty).\EN
\item If $\bar{x}$ is uniformly asymptotically stable then there exists a class $\mathcal{K}\mathcal{L}$ function $\beta$ and a neighborhood $\mathcal{N}_{\bar{x}}$, such that
\EQ&& d(\Phi_{f}(t,t_{0},x_{0}),\bar{x})\leq \beta(d(x_{0},\bar{x}),t-t_{0}),\hspace{.2cm}\nnum\\&& x_{0}\in \mathcal{N}_{\bar{x}}, t\in[t_{0},\infty).\EN
\end{itemize}
\end{theorem}
\begin{pf}
Let us consider a neighborhood $\mathcal{U}_{\bar{x}}\subset \exp_{\bar{x}}B_{i(\bar{x})}(0)$, where $i(\bar{x})$ is the injectivity radius at $\bar{x}\in M$  and $B_{i(\bar{x})}(0)\subset T_{\bar{x}}M$. Note that $i(\bar{x})>0$, see Proposition 2.1.10 in \cite{klin}. In order to prove the first assertion we note that the uniform Lyapunov stability of $\bar{x}$, implies that there exists $\mathcal{W}_{\bar{x}}\subset M$, such that $x_{0}\in \mathcal{W}_{\bar{x}}$ results in $\Phi_{f}(t,t_{0},x_{0})\in \mathcal{U}_{\bar{x}}$ for all $t\in [t_{0},\infty)$. Hence, $\mathcal{W}_{\bar{x}}\subseteq\mathcal{U}_{\bar{x}}\subset \exp_{\bar{x}}B_{i(\bar{x})}(0)$ and $\Phi_{f}(t,t_{0},x_{0})$ remains in a normal neighborhood of $\bar{x}$.

Lemma \ref{lpp} implies that $\exp_{\bar{x}}B_{r}(0)=B(\bar{x},r),$ provided $0< r\leq i(\bar{x})$. Hence, for any $\mathcal{U}_{\bar{x}}\doteq B(\bar{x},r),\hspace{.2cm}0<r< i(\bar{x})$, there exists $\mathcal{W}^{r}_{\bar{x}}\subseteq B(\bar{x},r)$, such that, $x_{0}\in \mathcal{W}^{r}_{\bar{x}}$ results in $\Phi_{f}(t,t_{0},x_{0})\in B(\bar{x},r),\hspace{.2cm}t\in[t_{0},\infty)$.
For any $0<r< i(\bar{x})$, define 
\EQ \label{n}&&\widehat{\mathcal{W}}_{\bar x}^r\doteq \mbox{int}\left(\{x\in M\hspace{0cm}|\hspace{0cm}\bigcup_{t\in[t_{0},\infty)}\Phi_{f}(t,t_{0},x)\in B(\bar{x},r)\}\right),\nnum\\&&\EN
where $\mbox{int}(A)$ is the interior set  of $A$. We note that $\mathcal{W}^{r}_{\bar{x}}$ is an open set and $\widehat{\mathcal{W}}_{\bar x}^r$ is the largest open set containing elements of $M$ for which the state trajectory $\Phi_{f}$ is contained in $B(\bar{x},r)$. Hence, $\mathcal{W}^{r}_{\bar{x}}\subseteq \widehat{\mathcal{W}}_{\bar x}^r\subseteq B(\bar{x},r)$. Since $\widehat{\mathcal{W}}_{\bar x}^r$ is an open set in 
        $M$ and the induced topology by the distance function $d$ is the same as the manifold topology (Theorem \ref{t2}), there always exists $l\in \mathds{R}_{>0}$, such that  $\exp_{\bar{x}}B_{l}(0)\subseteq \widehat{\mathcal{W}}_{\bar x}^r$.
Define
\EQ \label{del}\delta(r)\doteq \left\{ \begin{array}{ll}
         \max l\hspace{.2cm}|\exp_{\bar{x}}B_{l}(0)\subseteq \widehat{\mathcal{W}}_{\bar x}^r,\hspace{.7cm}r\leq i(\bar{x}),\\
        \max l\hspace{.2cm}|\exp_{\bar{x}}B_{l}(0)\subseteq \widehat{\mathcal{W}}_{\bar x}^{i(\bar{x})}\hspace{.5cm}i(\bar{x})< r.\end{array} \right.\nnum\\\EN
  Note that $l\in \mathds{R}_{\geq 0}$.       
Since our argument is local, without loss of generality, we assume $ i(\bar{x})<\infty$. Then for $r< i(\bar{x})$ we have $\overline{\widehat{\mathcal{W}}_{\bar x}^r}\subset \exp_{\bar{x}}\overline{B_{r}(0)}$, which together with the compactness of $\overline{B_{r}(0)}\subset T_{\bar{x}}M$ ($T_{\bar{x}}M$ is a finite dimensional vector space, $\overline{B_{r}(0)}$ is a closed and bounded set, and $\exp_{\bar{x}}$ is a local diffeomorphism), implies $0<\delta(r)<\infty$. Note that since $\exp$ is a local diffeomorphism then $\exp_{\bar{x}}\overline{B_{r}(0)}=\overline{\exp_{\bar{x}}B_{r}(0)}$. Now we show that $\delta(\cdot)$ is non-decreasing. Suppose $\delta(r)$ is strictly decreasing. Then, for $r_{1}<r_{2}< i(\bar{x})$ we have $\delta(r_{2})< \delta(r_{1})$. Denote the associated neighborhoods of $B(\bar{x},r_{1})$ and $B(\bar{x},r_{2})$ by $\widehat{\mathcal{W}}^{r_{1}}_{\bar{x}}$ and $\widehat{\mathcal{W}}^{r_{2}}_{\bar{x}}$ respectively, see \ref{n}. Then  $\delta(r_{2})< \delta(r_{1})$ implies that
\EQ \exists x_{0}\in \widehat{\mathcal{W}}^{r_{1}}_{\bar{x}},\hspace{.2cm}s.t.\hspace{.2cm}x_{0}\notin \widehat{\mathcal{W}}^{r_{2}}_{\bar{x}},\EN
where $B(\bar{x},r_{1})\subset B(\bar{x},r_{2})$. However, $x_{0}\in \widehat{\mathcal{W}}^{r_{1}}_{\bar{x}}$ results in $\Phi_{f}(t,t_{0},x_{0})\in B(\bar{x},r_{1})\subset B(\bar{x},r_{2}),\hspace{.2cm}t\in[t_{0},\infty)$, which contradicts $x_{0}\notin \widehat{\mathcal{W}}^{r_{2}}_{\bar{x}}$. Hence, $\delta(r_{1})\leq \delta(r_{2})$. 

Choose a $\zeta\in\mathcal{K}$ such that $\zeta(r)\leq \delta(r),\hspace{.2cm}r\in\mathds{R}_{\geq 0}$ (this is always possible since $\delta$ is non-decreasing), and $\zeta^{-1}:[0,\sup_{r\in[0,\infty)}\zeta(r))\mapsto \mathds{R}_{\geq 0}$ is a $\mathcal{K}$ class function. Note that $\zeta$ is bounded by $\delta$, hence, $\sup_{r\in[0,\infty)}\zeta(r)$ is bounded. Now choose  $\mathcal{N}_{\bar{x}}=\exp_{\bar{x}}B_{\sup_{r\in[0,\infty)}\zeta(r)}(0)\subset\exp_{\bar{x}}B_{\delta(i(\bar{x}))}(0)$. Then, $r\doteq \zeta^{-1}(d(x_{0},\bar{x})),\hspace{.2cm}x_{0}\in \mathcal{N}_{\bar{x}}$, implies
\EQ  d(x_{0},\bar{x})=\zeta(r)\leq \delta(r),\hspace{.2cm}x_{0}\in\mathcal{N}_{\bar{x}},\EN
and hence, by (\ref{n}) 
\EQ &&\forall x_{0}\in \mathcal{N}_{\bar{x}}\Rightarrow \Phi_{f}(t,t_{0},x_{0})\in B\left(\bar{x},r\right)=\nnum\\&& B\left(\bar{x},\zeta^{-1}(d(x_{0},\bar{x}))\right),t\in[t_{0},\infty).\EN
Hence, $d(\Phi_{f}(t,t_{0},x_{0}),\bar{x})\leq \zeta^{-1}(d(x_{0},\bar{x}))\doteq\alpha(d(x_{0},\bar{x})),\hspace{.2cm}\\t\in[t_{0},\infty)$, which proves the first statement. 
\\\\
The proof of the second assertion  follows from the proof given in \cite{Kha} by employing the function $\alpha$ constructed above.
\halmos
\end{pf}
\begin{remark}
Theorem \ref{lyp} characterizes the local behavior of state trajectories for uniformly stable/asymptotic stable dynamical systems on Riemannian manifolds. It has been shown that the Riemannian distance between the state trajectories and equilibriums are bounded above by positive continuous functions of the Riemannian distance function of the initial state. In the case of $M=\mathds{R}^{n}$, these properties recover the analogous stability properties of stable/ asymptotic stable dynamical systems on $\mathds{R}^{n}$, see \cite{Kha}, Chapter 4. 
\end{remark}
The following theorem gives the existence of Lyapunov functions and also characterizes their properties for locally asymptotically stable  systems evolving on Riemannian manifolds in normal neighborhoods of equilibriums of dynamical systems. In \cite{pai, simp} the Riemannian distance function is employed as a candidate to construct  Lyapunov functions for dynamical systems on Riemannian manifolds. For general discussions on the construction of Lyapunov functions on Riemannian manifolds and metric spaces see \cite{pai,nad,chic,simp,con,hur}.
\begin{theorem}
\label{lf1}
Let $\bar{x}$ be an equilibrium for the smooth dynamical system (\ref{peter}) on an open set $\mathcal{N}_{\bar{x}}\subset \mathcal{U}^{n}_{\bar{x}}$ ($\mathcal{U}^{n}_{\bar{x}}$ is a normal neighborhood around $\bar{x}$), such that there exists a $\mathcal{KL}$ function $\beta$, which satisfies
\EQ &&d(\Phi_{f}(t,t_{0},x_{0}),\bar{x})\leq \beta(d(x_{0},\bar{x}),t-t_{0}),\hspace{.2cm}\nnum\\&& x(t_{0})=x_{0}\in \mathcal{N}_{\bar{x}},t\in[t_{0},\infty).\EN
 Assume $||T_{x}f(\cdot,t)||$ is uniformly bounded with respect to $t$ on $\mathcal{N}_{\bar{x}}$, where $||.||$ is the norm of the bounded linear operator $Tf:TM\rightarrow TTM$ as per Definition \ref{d2}. Then, for some $\mathcal{U}_{\bar{x}}\subset \mathcal{U}^{n}_{\bar{x}}$, for all $x(t_{0})=x_{0}\in \mathcal{U}_{\bar{x}}$, there exist a Lyapunov candidate function $w:M\times \mathds{R}\rightarrow \mathds{R}_{\geq 0}$ and $\alpha_{1},\alpha_{2},\alpha_{3},\alpha_{4}\in \mathcal{K}$, such that for all $x\in \mathcal{U}_{\bar{x}}$ and $t\in[t_{0},\infty),$
\EQ \label{koonkoon1}(i):&&\hspace{.2cm}\alpha_{1}\left(d(x,\bar{x})\right)\leq w(x,t)\leq\alpha_{2}\left(d(x,\bar{x})\right),\nnum\\
 (ii):&&\hspace{.2cm}\mathfrak{L}_{f(x,t)}w\leq -\alpha_{3}\left(d(x,\bar{x})\right),\nnum\\ 
 (iii):&&\hspace{.2cm}||T_{x}w||\leq \alpha_{4}\left(d(x,\bar{x})\right),\EN
where $d(\cdot,\cdot)$ is the Riemannian metric, $\mathfrak{L}$ is the Lie derivative and $Tw:TM\rightarrow T\mathds{R}\simeq \mathds{R}\times \mathds{R}$ is the pushforward of $w$.
\end{theorem}

\begin{pf}
By employing Lemma \ref{eun}, consider $\exp_{\bar{x}}B_{\epsilon}(0)\subset M,\hspace{.2cm}0<\epsilon$, such that $\exp_{\bar{x}}$ is a diffeomorphism onto its image, then  $\exp_{\bar{x}}$ is invertible and the inverse map is denoted by $\exp^{-1}_{\bar{x}}:M\rightarrow T_{\bar{x}}M$. By Theorem \ref{t2} the induced topology of the distance function $d$ is the same as the original topology of $M$ and by Lemma \ref{lpp} the metric balls and geodesic balls are identical. Hence, without loss of generality, we assume $\mathcal{N}_{\bar{x}}=\exp_{\bar{x}}B_{\epsilon}(0)$, where $\exp_{\bar{x}}$ is a diffeomorphism onto $\exp_{\bar{x}}B_{\max\{\epsilon,\beta(\epsilon,0)\}}(0)$.

 Since $\exp_{\bar{x}}$ is a diffeomorphism onto $\exp_{\bar{x}}B_{\epsilon}(0)$, then for any $x\in \exp_{\bar{x}}B_{\epsilon}(0)$, there exists $z\in T_{\bar{x}}M$ such that $x=\exp_{\bar{x}}z$, or equivalently $z=\exp^{-1}_{\bar{x}}x$. Let us call the operator $\exp^{-1}_{\bar{x}}$ the \textit{geodesic lift}. The time variation of $z$, as long as $x$ stays in $\exp_{\bar{x}}B_{\epsilon}(0)$, is given by
\EQ \label{log}\dot{z}(t)&=&T_{x}\exp^{-1}_{\bar{x}}\left(f(x,t)\right)=T_{\exp_{\bar{x}}z}\exp^{-1}_{\bar{x}}\left(f(\exp_{\bar{x}}z,t)\right)\nnum\\&\doteq&\hat{f}(z,t),\nnum\\\EN
where $\dot{z}(t)\in T_{z(t)}T_{\bar{x}}M\simeq T_{\bar{x}}M$.
We note that the equilibrium $\bar{x}$ of $f(x,t)$ changes to $z=0\in T_{\bar{x}}M$ for the dynamical equations in $z$ coordinates. In the case $M=\mathds{R}^{n}$, we have
\EQ x=\exp_{\bar{x}}z=\bar{x}+z\in \mathds{R}^{n}.\EN 

For any $x(t_{0})\in\exp_{\bar{x}}B_{\epsilon}(0)$, we have $x(t_{0})=\exp_{\bar{x}}z(t_{0})$ for some $z(t_{0})\in B_{\epsilon}(0)$. Now let us consider the geodesic curve $\gamma:[0,1]\rightarrow M,\hspace{.2cm}\gamma(\tau)\doteq \exp_{\bar{x}}\tau z(t_{0})$. Employing the results of \cite{Lee3}, Proposition 5.11, in the normal coordinates of $\bar{x},$ we have
\EQ \gamma(\tau)=(\tau z_{1}(t_{0}),...,\tau z_{n}(t_{0})),\EN
where $d(x(t_{0}),\bar{x})=\left(\sum^{n}_{i=1}z^{2}_{i}(t_{0})\right)^{\frac{1}{2}}=||z(t_{0})||_{e}=||z(t_{0})||_{g}$. The last equality is due to the fact that in normal coordinates of $\bar{x}$, the Riemannian metric is given by 
\EQ g\left(\frac{\partial}{\partial x_{i}},\frac{\partial}{\partial x_{j}}\right)=\delta_{ij}+O(r^{2}),\EN
where $r$ is the distance and $\delta_{ij}$ is the Kronecker delta, see \cite{Pet}, Chapter 5. Hence, $g_{\bar{x}}(\frac{\partial}{\partial x_{i}},\frac{\partial}{\partial x_{j}})=\delta_{ij}$ and $||z(t_{0})||_{e}=||z(t_{0})||_{g}$.
Therefore, we have
\EQ &&||z(t)||_{g}\leq \beta(||z(t_{0})||_{g},t-t_{0}),\hspace{.2cm}\nnum\\&&z(t_{0})=z_{0}\in B_{\epsilon}(0)\subset T_{\bar{x}}M.\EN 
The uniform boundedness of $T_{x}f(\cdot,t)$ with respect to $t$ together with (\ref{log}) and smoothness of $\exp^{-1}$ imply that $\frac{\partial \hat{f}}{\partial z}$ is uniformly bounded on $B_{\epsilon}(0)\in T_{\bar{x}}M$. Hence, we can apply Theorem 4.16 of \cite{Kha} to demonstrate  the existence of a Lyapunov function $v:T_{\bar x}M \times\mathds{R}\rightarrow\mathds{R}$, satisfying
\EQ \label{nanat}(i):&&\hspace{.2cm}\alpha_{1}(||z||_{g})\leq v(z,t)\leq\alpha_{2}(||z||_{g}),\nnum\\
 (ii):&&\hspace{.2cm}\mathfrak{L}_{\hat{f}(z,t)}v\leq -\alpha_{3}(||z||_{g}),\nnum\\ (iii):&&\hspace{.2cm}|T_{z}v(\frac{\partial}{\partial z})|\leq \alpha_{4}(||z||_{g}),\EN 
where $ z\in T_{\bar{x}}M,t\in[t_{0},\infty)$.
 Since $\exp_{\bar{x}}$ is a local diffeomorphism by Lemma \ref{eun}, for $x\in\mathcal{N}_{\bar{x}}$, we have $x=\exp_{\bar{x}}\circ \exp^{-1}_{\bar{x}}x$. Hence, by \cite{Lee2}, Lemma 3.5
 \EQ \mbox{Id}&&=T_{x}\left(\exp_{\bar{x}}\circ\exp^{-1}_{\bar{x}}\right)=T_{\exp^{-1}_{\bar{x}}x}\exp_{\bar{x}}\circ T_{x}\exp^{-1}_{\bar{x}},\EN
 where $\mbox{Id}$ is the identity map and $T$ is the pushforward as per Definition \ref{d2}.
 This shows
 \EQ\label{kirkirkir} T_{x}\exp^{-1}_{\bar{x}}= \left(T_{\exp^{-1}_{\bar{x}}x}\exp_{\bar{x}}\right)^{-1}.\EN
 The Lie derivative of $v$ with respect to $\hat{f}$ is locally given by (\ref{lie}) as follows
 \EQ \mathfrak{L}_{\hat{f}(z,t)}v=dv(\frac{\partial}{\partial t},\hat{f}(z,t))=d_{t}v(\frac{\partial}{\partial t})+d_{z}v(\hat{f}(z,t)).\nnum\\\EN
  Since $v$ is a scalar-valued function then $dv(\frac{\partial}{\partial t},\hat{f}(z,t))=Tv(\frac{\partial}{\partial t},\hat{f}(z,t))=T_{t}v(\frac{\partial}{\partial t})+T_{z}(\hat{f}(z,t))$. Employing $\exp_{\bar{x}}$, we define the following function on $M$:
  \EQ \label{ch}\hat{v}(x,t)\doteq v(\exp^{-1}_{\bar{x}}x,t),\hspace{.2cm}x\in \exp_{\bar{x}}B_{\epsilon}(0).\EN
  
  Then the Lie derivative of $\hat{v}$ along $f$ at state $x$ and time $t$ is 
  \EQ \mathfrak{L}_{f(x,t)}\hat{v}&&=d_{t}\hat{v}(\frac{\partial}{\partial t})+T_{x}\hat{v}(f(x,t))\nnum\\&&=d_{t}\hat{v}(\frac{\partial}{\partial t})+T_{z}v\left(T_{x}\exp^{-1}_{\bar{x}}\circ T_{z}\exp_{\bar{x}}(\hat{f}(z,t))\right)\nnum\\&&=d_{t}\hat{v}(\frac{\partial}{\partial t})+T_{z}v\left(T_{x}\exp^{-1}_{\bar{x}}\circ \right.\nnum\\&&\left.T_{\exp^{-1}_{\bar{x}}x}\exp_{\bar{x}}(\hat{f}(z,t))\right)\nnum\\&&=d_{t}\hat{v}(\frac{\partial}{\partial t})+T_{z}v\left(\hat{f}(z,t)\right)\hspace{.25cm}\mbox{by employing (\ref{kirkirkir})}\nnum\\&&=d_{t}v(\frac{\partial}{\partial t})+T_{z}v\left(\hat{f}(z,t)\right)\hspace{.25cm}\mbox{by employing (\ref{ch})}\nnum\\&&=\mathfrak{L}_{\hat{f}(z,t)}v.\EN
  The same argument applies to $T_{x}\hat{v}(\frac{\partial}{\partial x})$ and shows that 
  \EQ T_{x}\hat{v}(\frac{\partial}{\partial x})=T_{z}v(\frac{\partial}{\partial z}).\EN
  
	As shown before we have $d(x(t),\bar{x})=||z(t)||_{g}$, hence, by (\ref{nanat}), $\hat{v}$ locally satisfies (\ref{koonkoon1}).
	
  Since the function constructed above is defined locally, it remains to extend the domain of its definition to $M$. 
  For $\delta\in (0,\epsilon)$, compactness of $\overline{B}_{\delta}(0)\subset T_{\bar{x}}M$ and smoothness of $\exp_{\bar{x}}$ together imply that $\exp_{\bar{x}}\overline{B}_{\delta}(0)\subset \exp_{\bar{x}}B_{\epsilon}(0)$ is a compact set in $M$. Choose a \textit{bump function} $\psi\in C^{\infty}(M)$, such that $\psi\equiv 1$ on $\exp_{\bar{x}}\overline{B}_{\delta}(0)$ and $\mbox{supp}\psi\subset \exp_{\bar{x}}B_{\epsilon}(0)$, where $\mbox{supp}\psi\doteq\{x\in M\hspace{.2cm}s.t.\hspace{.2cm}\psi(x)\ne 0\}$, for the definition of bump functions see \cite{Lee2}. As shown in \cite{Lee2}, Proposition 2.26, such bump functions always exist. Hence, we consider $\mathcal{U}_{\bar{x}}\doteq \exp_{\bar{x}}B_{\delta}(0)$ and $w\doteq \psi\times \hat{v}:M\times \mathds{R}\rightarrow \mathds{R}$. The Lie derivative of $w$ is given by
  \EQ \mathfrak{L}_{f(x,t)}w=\mathfrak{L}_{f(x,t)}\psi\cdot\hat{v}=\psi\mathfrak{L}_{f(x,t)}\hat{v}+\hat{v}\mathfrak{L}_{f(x,t)}\psi,\EN
  where on $\mathcal{U}_{\bar{x}}$ we have
  \EQ \mathfrak{L}_{f(x,t)}w=\mathfrak{L}_{f(x,t)}\hat{v}.\EN
  Same argument shows that on $\mathcal{U}_{\bar{x}}$, $T_{x}w(\frac{\partial}{\partial x})=T_{x}v(\frac{\partial}{\partial x})$, which completes the proof for the Lyapunov function $w$.
\halmos\end{pf}
Note that properties (ii) and (iii) are essential to obtain the robustness results for perturbed dynamical systems, see \cite{Kha}, Chapters 9,10,11. 
 The following theorem strengthens the hypotheses of Theorem \ref{lyp} to local exponential stability and derives the local properties of Lyapunov functions in a normal neighborhood of equilibriums. 
 
\begin{theorem}
\label{lffff}
Let $\bar{x}$ be a uniformly exponentially stable equilibrium of the dynamical system (\ref{peter}) on $\mathcal{N}_{\bar {x}}\subset \mathcal{U}^{n}_{\bar{x}}$ ($\mathcal{U}^{n}_{\bar{x}}$ is a normal neighborhood around $\bar{x}$), where $\mathcal{N}_{\bar {x}}$ denotes a subset of a 
normal neighborhood on an $n$ dimensional Riemannian manifold $(M,g)$. Assume  $||T_{x}f(\cdot,t)||$ is uniformly bounded, where $||.||$ is the norm of the linear operator $Tf:TM\rightarrow TTM$. Then, for some $\mathcal{U}_{\bar{x}}\subset \mathcal{U}^{n}_{\bar{x}}$, for all $x(t_{0})\in \mathcal{U}_{\bar{x}}$, there exist a Lyapunov function $v:M\times \mathds{R}\rightarrow \mathds{R}_{\geq 0}$ and $\lambda_{1},\lambda_{2},\lambda_{3},\lambda_{4}\in \mathds{R}_{>0}$, such that for all $x\in \mathcal{U}_{\bar{x}}$
\EQ \label{kirkir}(i):&&\hspace{.2cm}\lambda_{1}d^{2}(x,\bar{x})\leq v(x,t)\leq\lambda_{2} d^{2}(x,\bar{x}),\nnum\\
 (ii):&&\hspace{.2cm}\mathfrak{L}_{f(x,t)}v\leq -\lambda_{3}d^{2}(x,\bar{x}),\nnum\\ (iii):&&\hspace{.2cm}||T_{x}v||\leq \lambda_{4}d(x,\bar{x}).\EN\end{theorem}
\begin{pf}
Following the proof of Theorem \ref{lf1}, we employ the geodesic lift operator $z=\exp_{\bar{x}}^{-1}x$ in a normal neighborhood of $\bar{x}$. Hence, we obtain the local exponential stability of $0\in T_{\bar{x}}M$, for the dynamical system $\dot{z}(t)=\hat{f}(z,t)=T_{\exp_{\bar{x}}z}\exp_{\bar{x}}^{-1}\left(f(\exp_{\bar{x}}z,t)\right)$ as per the proof of Theorem \ref{lf1}. 
The rest of the proof parallels the proof of Theorem \ref{lf1} and the results of \cite{Kha}, Theorem 4.14. 
\halmos\end{pf}
We note that by employing the normal coordinates used in the proof of Theorem \ref{lf1}, we have $d(\Phi_{f}(t,t_{0},x(t_{0})),\bar{x})\leq Kd(x(t_{0}),\bar{x})\exp(-\lambda(t-t_{0}))$ implies \\$||z(t)||_{g}\leq K\exp(-\lambda(t-t_{0}))||z(t_{0})||_{g}$ which is required in the proof of Theorem 4.14 in \cite{Kha}.

%\subsection{Converse Lyapunov Theorem on Coordinate Charts of $M$}
The Lyapunov functions in Theorems \ref{lf1} and \ref{lffff} are constructed in a normal neighborhood of an equilibrium where $\exp_{x}$ is a local diffeomorphism. Hence, the properties derived in  Theorems \ref{lf1} and \ref{lffff} hold locally and the corresponding neighborhoods are restricted by the injectivity radius of the equilibrium. Depending on the geometric features of $M$, the injectivity radius of a particular point might be very small. In this section we construct  Lyapunov functions on a compact subset of a local chart of an equilibrium of a dynamical system on $M$ by scaling the Riemannian and Euclidean metrics. This is also a local method since we are restricted to work within a local coordinate system. However, in some cases, it may provide much larger neighborhood on which Theorems \ref{lf1} and \ref{lffff} hold.

    \begin{theorem}
\label{lf2}
Let $\bar{x}$ be an equilibrium for the dynamical system (\ref{peter}) on a coordinate chart $(U,\phi)$ of $\bar{x}$ as per Definition \ref{cc}, such that there exists a $\mathcal{KL}$ function $\beta$, which satisfies
\EQ \label{dr}d(\Phi_{f}(t,t_{0},x_{0}),\bar{x})\leq \beta(d(x_{0},\bar{x}),t-t_{0}),\hspace{.2cm}x(t_{0})=x_{0}\in U.\nnum\\\EN
 Assume $||T_{x}f(\cdot,t)||$ is uniformly bounded with respect to $t$ on $U$, where $||.||$ is the norm of the linear operator $Tf(x,t):TM\rightarrow TTM$. Then, for some $\mathcal{U}_{\bar{x}}\subset U$, for all $x(t_{0})=x_{0}\in \mathcal{U}_{\bar{x}}$, there exist a Lyapunov function $w:M\times \mathds{R}\rightarrow \mathds{R}_{\geq 0}$ and $\alpha_{1},\alpha_{2},\alpha_{3},\alpha_{4}\in \mathcal{K}$, such that for all $x\in \mathcal{U}_{\bar{x}}$
\EQ \label{koonkoon}(i):&&\hspace{.2cm}\alpha_{1}\left(d(x,\bar{x})\right)\leq w(x,t)\leq\alpha_{2}\left(d(x,\bar{x})\right),\nnum\\
 (ii):&&\hspace{.2cm}\mathfrak{L}_{f(x,t)}w\leq -\alpha_{3}\left(d(x,\bar{x})\right),\nnum\\ 
 (iii):&&\hspace{.2cm}||T_{x}w||\leq \alpha_{4}\left(d(x,\bar{x})\right).\EN
\end{theorem}
\begin{pf}
For the coordinate chart $(U,\phi)$, we have $\phi:M\rightarrow \mathds{R}^{n}$. By definition, $\phi$ is a homeomorphism to an open set in $\mathds{R}^{n}$, see \cite{Lee4}. Without loss of generality, we assume $\phi(\bar{x})=(0,...,0)$, otherwise we can consider the map $\eta(x)\doteq \phi(x)-\phi(\bar{x})$, where $\eta$ is also a homeomorphism by definition. Denote
\EQ &&\mathfrak{R}\doteq \max r,\, s.t. \hspace{.2cm}B_{e}(r,0)\subset \phi(U), \phi^{-1}\left(\overline{B_{e}(r,0)}\right)\subset U,\nnum\\&& r\in \mathds{R}_{>0},\EN
where $B_{e}(r,0)$ is the Euclidean ball of radius $r$. In $\mathds{R}^{n}$, $\overline{B_{e}(r,0)}$ is a compact set and since $\phi$ is a homeomorphism then, $\phi^{-1}\left(\overline{B_{e}(r,0)}\right)\subset M$ is a compact set. By (\ref{dr}), there exists $\mathcal{W}_{\bar{x}}\subset M$, such that for all $x_{0}\in \mathcal{W}_{\bar{x}} $,  $\Phi(t,t_{0},x_{0})\in \phi^{-1}\left(B_{e}(\mathfrak{R},0)\right)$. 

Replace the Riemannian metric $||\cdot||_{g}$ by the Euclidean metric $||\cdot||_{e}$ on $\phi^{-1}\left(\overline{B_{e}(\mathfrak{R},0)}\right)$. Since $\phi^{-1}\left(\overline{B_{e}(\mathfrak{R},0)}\right)$ is compact, there exists $c_{1},c_{2}\in \mathds{R}_{>0}$, such that \cite{Lee3}
\EQ\label{sc} c_{1}||X||_{g}\leq ||X||_{e}\leq c_{2}||X||_{g},\, \nnum\\X\in T_{x}M,x\in \phi^{-1}\left(\overline{B_{e}(\mathfrak{R},0)}\right).\EN
Since the state trajectory is contained in $\phi^{-1}\left(B_{e}(\mathfrak{R},0)\right)$, by replacing the Riemannian metric with the Euclidean one, the state trajectory will be bounded in $B_{e}(\mathfrak{R},0)$. By employing (\ref{sc}), the Euclidean distance function is bounded by the Riemannian one as follows.
Consider any piecewise smooth curve $\gamma:[a,b]\rightarrow 
M$ connecting $x\in \phi^{-1}(B_{e}(\mathfrak{R},0))$ and $\bar{x}$, such that $\gamma(a)=x$ and $\gamma(b)=\bar{x}$. If $\gamma$ belongs to $\phi^{-1}\left(\overline{B_{e}(\mathfrak{R},0)}\right)\subset M$, then
\EQ ||x-\bar{x}||_{e}\leq \int^{b}_{a}||\dot{\gamma}(s)||_{e}ds\leq c_{2}\int^{b}_{a}||\dot{\gamma}(s)||_{g}ds,\EN
where $||x-\bar{x}||_{e}$ is the Euclidean distance between $x$ and $\bar{x}$. In case $\gamma$ does not entirely belong to $\phi^{-1}\left(\overline{B_{e}(\mathfrak{R},0)}\right)$, then there exists a time $t\in [a,b]$, such that $\gamma(s)\in \phi^{-1}\left(S_{e}(\mathfrak{R},0)\right), s\in[a,t]$ and $||x-\gamma(t)||_{e}=\mathfrak{R}$, where $S_{e}(\mathfrak{R},0)=\{x\,|\,||x||_{e}=\mathfrak{R}\}$. Hence, since $x\in \phi^{-1}\left(B_{e}(\mathfrak{R},0)\right)$, we have
\EQ ||x-\bar{x}||_{e}&&\leq\mathfrak{R}\leq \int^{t}_{a}||\dot{\gamma}(s)||_{e}ds\leq c_{2}\int^{t}_{a}||\dot{\gamma}(s)||_{g}ds\nnum\\
&&\leq c_{2}\int^{b}_{a}||\dot{\gamma}(s)||_{g}ds.\EN
Therefore, for any piecewise smooth $\gamma$, $||x-\bar{x}||_{e}\leq c_{2}\int^{b}_{a}||\dot{\gamma}(s)||_{g}ds$. Taking the infimum over all $\gamma$, (\ref{len}) implies that
\EQ ||x-\bar{x}||_{e}\leq c_{2}d(x,\bar{x}).\EN
As $B_{e}(\mathfrak{R},0)$ is a convex set, the line connecting $x$ and $\bar{x}$ is entirely in $\phi^{-1}\left(B_{e}(\mathfrak{R},0)\right)$. Hence,
\EQ c_{1}d(x,\bar{x})\leq c_{1}\int^{b}_{a}||\dot{\gamma}(s)||_{g}ds&&\leq \int^{b}_{a}||\dot{\gamma}(s)||_{e}ds\nnum\\&&=||x-\bar{x}||_{e}.\EN
Therefore, for $x=\Phi_{f}(t,t_{0},x_{0})$, we have
\EQ \frac{1}{c_{2}}||x-\bar{x}||_{e}&&\leq d(\Phi_{f}(t,t_{0},x_{0}),\bar{x})\leq \beta(d(x_{0},\bar{x}),t-t_{0})\nnum\\&&\leq \beta(\frac{1}{c_{1}}||x-\bar{x}||_{e},t-t_{0}).\EN
Hence, $||x-\bar{x}||_{e}\leq c_{2}\beta(\frac{1}{c_{1}}||x-\bar{x}||_{e},t-t_{0})\doteq\hat{\beta}(||x-\bar{x}||_{e},t-t_{0})$.

The Euclidan induced norm of $T_{x}f(x,t)$ is defined by
\EQ ||T_{x}f(x,t)||_{e}&&=\sup_{X\in T_{x}M, X\ne 0}\frac{||T_{x}f(x,t)(X)||_{e}}{||X||_{e}}\nnum\\&&\leq \sup_{X\in T_{x}M, X\ne 0}\frac{c_{2}||T_{x}f(x,t)(X)||_{g}}{c_{1}||X||_{g}}\nnum\\&&\leq \frac{c_{1}}{c_{2}}||T_{x}f(x,t)||_{g}.\EN
Hence, boundedness of $||T_{x}f(\cdot,t)||_{g}$ implies the boundedness of $||T_{x}f(\cdot,t)||_{e}$. We apply the results of \cite{Kha}, Theorem 4.16 to the dynamical system evolving on $M$, where $||\cdot||_{g}$ is replaced by $||\cdot||_{e}$. Therefore, there exist a Lyapunov function $v$ and functions $\alpha_{1},\alpha_{2},\alpha_{3},\alpha_{4}\in\mathcal{K}$, such that 
\EQ (i):&&\hspace{.2cm}\alpha_{1}(||x||_{e})\leq v(x,t)\leq\alpha_{2}(||x||_{e}),\nnum\\
 (ii):&&\hspace{.2cm}\mathfrak{L}_{f(x,t)}v\leq -\alpha_{3}(||x||_{e}),\nnum\\ (iii):&&\hspace{.2cm}||T_{x}v||_{e}\leq \alpha_{4}(||x||_{e}),\hspace{.2cm} x\in B_{e}(\mathfrak{R},0),\EN
 where $||x||_{e}=||x-\bar{x}||_{e}$.  As a result of the scaling the Riemannian and Euclidean norms, we have
 \EQ (i):&&\hspace{.2cm}\alpha_{1}(c_{1}d(x,\bar{x}))\leq v(x,t)\leq\alpha_{2}(c_{2}d(x,\bar{x})),\nnum\\
 (ii):&&\hspace{.2cm}\mathfrak{L}_{f(x,t)}v\leq -\alpha_{3}(c_{1}d(x,\bar{x})),\nnum\\ (iii):&&\hspace{.2cm}||T_{x}v||\leq\frac{c_{2}}{c_{1}} \alpha_{4}(c_{2}d(x,\bar{x})).\EN
 Following the last part of the proof of Theorem \ref{lf1}, the domain of the definition of $v$ can be extended to $M$, which completes the proof for $\mathcal{U}_{\bar{x}}=\mathcal{W}_{\bar{x}}$.
\halmos\end{pf}

\section {Stability of Perturbed Dynamical Systems}
\label{s4}
The properties of the constructed Lyapunov functions in Theorems \ref{lf1} and \ref{lffff} are employed to obtain the robust stability results for perturbed dynamical systems on Riemannian manifolds.  
Consider the following perturbed dynamical system on $(M,g)$.
\EQ \label{per}\dot{x}(t)=f(x,t)+h(x,t),f,h\in \mathfrak{X}(M\times \mathds{R}).\EN
The term $h$ can be considered as a perturbation of the nominal system $f$. As stated in \cite{Kha,mic,gou}, stability results for (\ref{per}) can be obtained based on technical assumptions on the stability of the nominal system $f$ and boundedness of $h$. %The following Lemmas gives the existence of Lyapunov functions in a normal neighborhood of an asymptotically  stable system.
  The following theorem gives the stability of (\ref{per}), where the nominal system is  locally uniformly asymptotically stable.
  \begin{theorem}
  \label{t7}
  Let $\bar{x}$ be an equilibrium of dynamical system (\ref{peter}), which is locally uniformly asymptotically stable on a normal neighborhood $\mathcal{N}_{\bar{x}}$. Assume the perturbed dynamical system (\ref{per}) is complete and the  Riemannian norm of the perturbation $h\in \mathfrak{X}(M\times \mathds{R})$ is bounded on $\mathcal{N}_{\bar{x}}$, i.e. $||h(x,t)||_{g}\leq \delta,x\in \mathcal{N}_{\bar{x}},t\in[t_{0},\infty)$. Then, for sufficiently small $\delta$, there exists a neighborhood $U_{\bar{x}}$ and a function $\rho\in\mathcal{K}$, such that 
  \EQ \limsup_{t\rightarrow \infty} d(\Phi_{f+h}(t,t_{0},x_{0}),\bar{x})\leq\rho(\delta),\hspace{.2cm}x_{0}\in U_{\bar{x}}.\EN  
  \end{theorem}
  \begin{pf}
  Following the proof of Theorem \ref{lf1}, there exists $\mathcal{U}_{\bar{x}}\subset\mathcal{N}_{\bar{x}}$, such that (\ref{koonkoon1}) holds for a Lyapunov function $w$. First we show that the neighborhood $\mathcal{U}_{\bar{x}}$ in Theorem \ref{lf1} can be shrunk, such that $\Phi_{f+h}(\cdot,t_{0},x_{0})\in \mathcal{U}_{\bar{x}}$ provided $x_{0}\in \mathcal{U}_{\bar{x}}$. By Lemma \ref{lc} there exists $\mathcal{N}_{b,t_{0}}(\bar{x})$ and $\alpha_{3}\in\mathcal{K}$, such that
  \EQ &&\mathfrak{L}_{f+h}w=\mathfrak{L}_{f}w+\mathfrak{L}_{h}w\leq -\alpha_{3}(d(x,\bar{x}))+\mathfrak{L}_{h}w,\nnum\\&&\hspace{.2cm}x\in int(\mathcal{N}_{b,t_{0}}(\bar{x})).\EN 
  By the Shrinking Lemma \cite{Lee4} there exists a precompact neighborhood $\mathcal{W}_{\bar{x}}$, such that, $\mathcal{W}_{\bar{x}}\subset int(\mathcal{N}_{b,t_{0}}(\bar{x}))\subset \mathcal{N}_{b,t_{0}}(\bar{x})$, see \cite{Lee4}. Hence, $M-\mathcal{W}_{\bar{x}}$ is a closed set and $ \mathcal{N}_{b,t_{0}}(\bar{x})\bigcap (M-\mathcal{W}_{\bar{x}})$ is a compact set (closed subsets of compact sets are compact). The continuity of $\alpha_{3}$ and $d(\cdot,\bar{x})$ together with the compactness of $ \mathcal{N}_{b,t_{0}}(\bar{x})\bigcap (M-\mathcal{W}_{\bar{x}})$ imply the existence of the following parameter $\mathfrak{M}$,
  \EQ \mathfrak{M}\doteq \sup_{x\in \mathcal{N}_{b,t_{0}}(\bar{x})\bigcap (M-\mathcal{W}_{\bar{x}})}-\alpha_{3}(d(x,\bar{x}))< 0.\EN
  Note that $\alpha_{3}\in \mathcal{K}$, $x\in \mathcal{N}_{b,t_{0}}(\bar{x})\bigcap (M-\mathcal{W}_{\bar{x}})$ and since $\mathcal{W}_{\bar{x}}$ is a neighborhood of $\bar{x}$ then $d(x,\bar{x})>0,\hspace{.2cm}x\in \mathcal{N}_{b,t_{0}}(\bar{x})\bigcap (M-\mathcal{W}_{\bar{x}})$. Therefore, $\mathfrak{M}<0$. Using (\ref{lie}) implies that $\mathfrak{L}_{h}w=dw(h)\leq ||dw||\cdot||h||_{g}\leq \delta ||dw||$, where $||dw||$ is the induced norm of the linear operator $dw:TM\rightarrow \mathds{R}$. The smoothness of $w$ and compactness of $\mathcal{N}_{b,t_{0}}(\bar{x})$ together imply $||dw||<\infty$. It is important to note that $||dw||$ is closely related to $||Tw||$ through the component of the Riemannian metric $g$. As is shown by Theorem \ref{lf1}, $||T_{x}w||\leq\alpha_{4}(d(x,\bar{x}))$. Hence, the smoothness of $M$ and compactness of $\mathcal{N}_{b,t_{0}}(\bar{x})$ imply that $||dw||<\infty$. Note that $||dw||$ is the norm of the linear operator $dw:T_{x}M\rightarrow \mathds{R}$. Hence, for sufficiently small $\delta$, we have $\mathfrak{L}_{f+h}w<0, x\in \mathcal{N}_{b,t_{0}}(\bar{x})\bigcap(M-\mathcal{W}_{\bar{x}})$. Therefore, the state trajectory $\Phi_{f+h}(\cdot,t_{0},x_{0})$ stays in $\mathcal{U}_{\bar{x}}$ for all $x_{0}\in int(\mathcal{N}_{b,t_{0}}(\bar{x}))$.
  
   Without loss of generality,  assume $\mathcal{U}_{\bar{x}}=\exp_{\bar{x}}B_{r_{2}}(0), r_{2}<i(\bar{x})$. Then, by the results of Theorem \ref{lf1}, the variation of $w$ along $f+h$ is then given by
  \EQ  \mathfrak{L}_{f+h}w&&=\mathfrak{L}_{f}w+\mathfrak{L}_{h}w\leq -\alpha_{3}(d(x,\bar{x}))+\mathfrak{L}_{h}w\nnum\\&&=-\alpha_{3}(d(x,\bar{x}))+dw(h(x,t))\nnum\\&&=-\alpha_{3}(d(x,\bar{x}))+T_{x}w(h(x,t)) \nnum\\&&\leq -\alpha_{3}(d(x,\bar{x}))+||T_{x}w||\cdot||h(x,t)||_{g}\nnum\\&&\leq -\alpha_{3}(d(x,\bar{x}))+\delta \alpha_{4}(d(x,\bar{x}))\nnum\\&&\leq -(1-\theta)\alpha_{3}(d(x,\bar{x}))-\theta\alpha_{3}(d(x,\bar{x}))\nnum\\&&+\delta \alpha_{4}(d(x,\bar{x}))\leq -(1-\theta)\alpha_{3}(d(x,\bar{x})),\hspace{.2cm}\nnum\\&&\mbox{if}\hspace{.2cm} \alpha^{-1}_{3}(\frac{\delta\alpha_{4}(r_{1})}{\theta})\leq d(x,\bar{x})\leq r_{2},\EN
  where $r_{1}<r_{2}$, $0<\theta<1$ and $\delta\leq \theta\frac{\alpha_{3}(r_{2})}{\alpha_{4}(r_{1})}$. 
  
  Define $\eta\doteq \alpha_{2}(\alpha^{-1}_{3}(\frac{\delta\alpha_{4}(r_{1})}{\theta}))$, then $\{x\in M\,|\,d(x,\bar{x})\leq \alpha^{-1}_{3}(\frac{\delta\alpha_{4}(r_{1})}{\theta})\}\subset\mathcal{N}_{t,\eta}=\{x\in M |\,w(x,t)\leq \eta\}\subset \{x\in M|\,\alpha_{1}(d(x,\bar{x}))\leq \eta\}$.
  Hence, solutions initialized in $\{x\in M\,|\,d(x,\bar{x})\leq \alpha^{-1}_{3}(\frac{\delta\alpha_{4}(r_{1})}{\theta})\}$ remain in $\{x\in M|\,\alpha_{1}(d(x,\bar{x}))\leq \eta\}$ since $\dot{w}<0$ for $x\in \mathcal{N}_{t,\eta}-\{x\in M\,|\,d(x,\bar{x})\leq \alpha^{-1}_{3}(\frac{\delta\alpha_{4}(r_{1})}{\theta})\}$. This proves 
  \EQ \limsup_{t\rightarrow \infty} d(\Phi_{f+h}(t,t_{0},x_{0}),\bar{x})&\leq& \alpha_{1}^{-1}\left(\alpha^{-1}_{3}\left(\frac{\delta\alpha_{4}(r_{1})}{\theta}\right)\right)\nnum\\&\doteq&\rho(\delta),\EN
  for any $x_{0}\in U_{\bar{x}}\doteq \{x\in M\,|\,d(x,\bar{x})< \alpha^{-1}_{3}(\frac{\delta\alpha_{4}(r_{1})}{\theta})\}\bigcap\\ int(\mathcal{N}_{b,t_{0}}(\hat{x}))$. 
 \halmos
\end{pf}  

In the following theorem we strengthen the uniform asymptotic stability to the uniform exponential stability for the nominal system $\dot{x}=f(x,t)$. It will be shown that the state trajectory of the  perturbed system stays close to the equilibrium of the nominal system when some specific conditions are satisfied. 
\begin{figure}
\begin{center}
%\hspace*{0cm}\includegraphics[scale=.3]{lyp}
 %\vspace*{-.5cm}\caption{Uniformly normal neighborhood of $\bar{x}$ and a sublevel set of $v$}
     \label{ff2}
      \end{center}
   \end{figure}   
\begin{theorem}
Let $\bar{x}$ be an equilibrium of (\ref{peter}), which is locally uniformly exponentially stable on a normal neighborhood  $\mathcal{U}_{\bar{x}}$. Assume the nominal and perturbed dynamical systems are both complete and the  Riemannian norm of the perturbation $h\in \mathfrak{X}(M\times \mathds{R})$ is bounded on $\mathcal{U}_{\bar{x}}$, i.e. $||h(x,t)||_{g}\leq \delta,x\in \mathcal{U}_{\bar{x}},t\in[t_{0},\infty)$. Also assume $||f||_{g}$ and $||Tf||$ are uniformly bounded with respect to $t$ on compact subsets of $M$, where $Tf:TM\rightarrow TTM$ as per Definition \ref{d2}.   
 Then, for sufficiently small $\delta$, there exists positive constants $\zeta,\gamma$ and $k$, such that 
  \EQ &&d(\Phi_{f+h}(t,t_{0},x_{0}),\bar{x})\leq k\exp(-\gamma(t-t_{0}))d(x_{0},\bar{x})+\zeta\delta.\nnum\\\EN 
\end{theorem}
\begin{pf}
By Theorem \ref{lffff} there exists a Lyapunov candidate function $v$ which satisfies (\ref{kirkir}). Hence,  following the proof of Theorem \ref{t7}, there exists a connected compact sublevel set of $v$, such that $\mathcal{N}_{b,t_{0}}(\bar{x})\subset \mathcal{U}_{\bar{x}}$, where $\Phi_{f+h}(t,t_{0},\hat{x}_{0})\in \mathcal{N}_{b,t_{0}}(\bar{x}), \hspace{.2cm}\hat{x}_{0}\in int(\mathcal{N}_{b,t_{0}}(\bar{x})),t\in[t_{0},\infty)$. Since $int(\mathcal{N}_{b,t_{0}}(\bar{x}))$ is an open set, for a given $x_{0}\in int(\mathcal{N}_{b,t_{0}}(\bar{x}))$, we can choose $\hat{x}_{0}$ sufficiently close to $x_{0}$, such that $\hat{x}_{0}\in int(\mathcal{N}_{b,t_{0}}(\bar{x}))$. 

%By the results of Lemma \ref{un}, any normal neighborhood of $\bar{x}$ contains a uniformly normal neighborhood $U_{\bar{x}}$ and by Lemma \ref{lc} any neighborhood of $\bar{x}$ contains a compact connected sublevel set $\mathcal{N}_{b}(\bar{x})$. Hence, without loss of generality, we assume $\mathcal{N}_{b}(\bar{x})$ is a subset of $U_{\bar{x}}$. Since $U_{\bar{x}}$ is a uniformly normal neighborhood of $\bar{x}$, then for $x_{0},\hat{x}_{0}\in int(\mathcal{N}_{b}(\bar{x})$, there exists $\delta>0$, such that $x_{0}\in \exp_{\hat{x}_{0}}B_{\delta}(0)$ and $\hat{x}_{0}\in \exp_{x_{0}}B_{\delta}(0)$, see Figure \ref{ff2}. Note that $B_{\delta}(0)$ is a subset of $T_{x_{0}}M$ and $T_{\hat{x}_{0}}M$ with respect to metrics $g_{x_{0}}$ and $g_{\hat{x}_{0}}$. Since $\Phi_{f+h}(t,t_{0},\hat{x}_{0})\in U_{\bar{x}},\,t\in[t_{0},\infty)$, then 

%\EQ \Phi_{f+h}(t,t_{0},\hat{x}_{0})=\exp_{\Phi_{f}(t,t_{0},x_{0})}z(t),\,t\in[t_{0},\infty),\EN
%where $z(t)\in T_{\Phi_{f}(t,t_{0},x_{0})}M$. Since $\exp$ is a local diffeomorphism then $z(t)=\exp^{-1}_{\Phi_{f}(t,t_{0},x_{0})}\Phi_{f+h}(t,t_{0},\hat{x}_{0})$. For simplicity in notation, denote $x(t)=\Phi_{f}(t,t_{0},x_{0})$ and $y(t)=\Phi_{f+h}(t,t_{0},\hat{x}_{0})$, where $x(t_{0})=x_{0}$ and $y(t_{0})=\hat{x}_{0}$. The time variation of $z=\exp^{-1}_{x}y$ is then given by
%\EQ \dot{z}(t)&&=T_{y(t)}\exp^{-1}_{x(t)}(\dot{y}(t))+T_{x(t)}\exp^{-1} y(t)\left(\dot{x}(t)\right)\nnum\\&&=T_{y(t)}\exp^{-1}_{x(t)}(f(y(t))+h(y(t)))+\nnum\\&&T_{x(t)}\exp^{-1} y(t)\left(f(x(t))\right).\EN
By employing the results of Theorem \ref{lffff}, the variation of $v$ along $f+h$ is then given by
  \EQ  \mathfrak{L}_{f+h}v&&=\mathfrak{L}_{f}v+\mathfrak{L}_{h}v\leq -\lambda_{3}d^{2}(x,\bar{x})+\mathfrak{L}_{h}v\nnum\\&&\leq -\lambda_{3}d^{2}(x,\bar{x})+\mathfrak{L}_{h}v=-\lambda_{3}d^{2}(x,\bar{x})+dv(h(x,t))\nnum\\&&=-\lambda_{3}d^{2}(x,\bar{x})+T_{x}v(h(x,t)) \nnum\\&&\leq -\lambda_{3}d^{2}(x,\bar{x})+||T_{x}v||\cdot||h(x,t)||_{g}\nnum\\&&\leq -\lambda_{3}d^{2}(x,\bar{x})+\delta \lambda_{4}d(x,\bar{x}).\EN
  
  Hence,
  \EQ \mathfrak{L}_{f+h}v=\dot{v}\leq -\frac{\lambda_{3}}{\lambda_{2}}v+\delta\lambda_{4}\sqrt{\frac{v}{\lambda_{1}}}.\EN
  Following the comparison method presented in \cite[Section 9.3]{Kha}, we have
  \EQ \sqrt{v(x,t)}&&\leq \sqrt{v(x_{0},t_{0})}\exp(-\frac{\lambda_{3}}{2\lambda_{2}}(t-t_{0}))\nnum\\&&+\delta\frac{\lambda_{4}\lambda_{2}}{\lambda_{3}\lambda_{1}}\Big[1-\exp(-\frac{\lambda_{3}}{2\lambda_{2}}(t-t_{0}))\Big].\EN
  Therefore, 
  \EQ d(\Phi_{f+h}(t,t_{0},x_{0}),\bar{x})&&\leq \sqrt{\frac{\lambda_{2}}{\lambda_{1}}}\exp(-\frac{\lambda_{3}}{2\lambda_{2}}(t-t_{0}))d(x_{0},\bar{x})\nnum\\&&+\delta\frac{\lambda_{4}\lambda_{2}}{\lambda_{3}\lambda_{1}}\Big[1-\exp(-\frac{\lambda_{3}}{2\lambda_{2}}(t-t_{0}))\Big]\nnum\\&&\leq \sqrt{\frac{\lambda_{2}}{\lambda_{1}}}\exp(-\frac{\lambda_{3}}{2\lambda_{2}}(t-t_{0}))d(x_{0},\bar{x})\nnum\\&&+\delta\frac{\lambda_{4}\lambda_{2}}{\lambda_{3}\lambda_{1}},\hspace{.2cm} \mbox{if}\, \delta\leq \frac{\lambda_{3}}{\lambda_{4}}\sqrt{\frac{\lambda_{1}}{\lambda_{2}}}d(x_{0},\bar{x}),\nnum\\\EN
  which completes the proof for $k\doteq \sqrt{\frac{\lambda_{2}}{\lambda_{1}}}, \gamma\doteq \frac{\lambda_{3}}{2\lambda_{2}}$ and $\zeta \doteq \frac{\lambda_{4}\lambda_{2}}{\lambda_{3}\lambda_{1}}$.
\halmos\end{pf}
\section{Conclusion}
\label{s5}
In this paper we have presented the stability results for dynamical systems evolving on Riemannian manifolds. We have obtained converse Lyapunov theorems for nonlinear dynamical systems defined on smooth connected Riemannian manifolds and have characterized properties of  Lyapunov functions with respect to the Riemannian distance function. The results are given by using the geometrical concepts such as normal neighborhoods, injectivity radius and bump functions on Riemannian manifolds.
\bibliographystyle{plain}        % Include this if you use bibtex 
\bibliography{HSCC1}           % and a bib file to produce the 
                             % bibliography (preferred). The
                                 % correct style is generated by
                                 % Elsevier at the time of printing.

%\begin{thebibliography}{99}     % Otherwise use the  
                                 % thebibliography environment.
                                 % Insert the full references here.
                                 % See a recent issue of Automatica 
                                 % for the style.
%  \bibitem[Heritage, 1992]{Heritage:92}
%     (1992) {\it The American Heritage. 
%     Dictionary of the American Language.}
%     Houghton Mifflin Company.
%  \bibitem[Able, 1956]{Abl:56}
%     B.~C.~Able (1956). Nucleic acid content of macroscope. 
%     {\it Nature 2}, 7--9. 
%  \bibitem[Able {\em et al.}, 1954]{AbTaRu:54}   
%     B.~C. Able, R.~A. Tagg, and M.~Rush (1954).
%     Enzyme-catalyzed cellular transanimations.
%     In A.~F.~Round, editor, 
%     {\it Advances in Enzymology Vol. 2} (125--247). 
%     New York, Academic Press.
%  \bibitem[R.~Keohane, 1958]{Keo:58}
%     R.~Keohane (1958).
%     {\it Power and Interdependence: 
%     World Politics in Transition.}
%     Boston, Little, Brown \& Co.
%  \bibitem[Powers, 1985]{Pow:85}
%     T.~Powers (1985).
%     Is there a way out?
%     {\it Harpers, June 1985}, 35--47.

%\end{thebibliography}
\appendix
%\section{A summary of Latin grammar}    % Each appendix must have a short title.
%\section{Some Latin vocabulary}         % Sections and subsections are supported  
                                     % in the appendices.
\end{document}